\documentclass[12pt, a4paper, dvips, dvipdf]{article}
\usepackage{amsmath,amssymb}

\begin{document}
\vspace{0cm}

\newenvironment{proof}[1][Proof]{\noindent {\it{#1.}}}{\ \rule{0.5em}{0.5em}}

\newtheorem{thm}{Theorem}[section]
\newtheorem{pro}[thm]{Proposition}
\newtheorem{cor}[thm]{Corollary}
\newtheorem{lem}[thm]{Lemma}
\newtheorem{defn}[thm]{Definition}
\newtheorem{exm}[thm]{Example}
\newtheorem{rem}[thm]{Remark}
\baselineskip=6.5mm

\title{\textbf{Strong geodesic convex functions of order m}}
\author{ Akhlad Iqbal$^*$ and Izhar Ahmad$^{**}$ \\
\\\textit{$^*$Department of Mathematics,Aligarh Muslim University,}
\\\textit{Aligarh-202002, India}
\\\textit{$^{**}$Department of Mathematics and Statistics, King Fahd University of}
\\\textit{Petroleum and Minerals, Dhahran 31261, Saudi Arabia}\\
\\\textit{emails:akhlad6star@gmail.com; drizhar@kfupm.edu.sa}}
\maketitle

 \noindent\textbf{Abstract.} Strong geodesic  convex function and strong monotone vector field of order $m$ on Riemannian manifolds have been established. A characterization of strong geodesic convex function of order $m$ for the continuously differentiable functions has been discussed. The relation between the solution of a new variational inequality problem and the strict minimizers of order $m$ for a multiobjective programming problem has also been established. \\

\noindent\textbf{2010 Mathematics Subject Classification: } 33B15; 26B25; 26D15; 26D10.

\noindent\textit{Keywords:}  Geodesic convex function, Monotone vector field, Variational inequality problem, Riemannian manifolds.\\

\noindent
\section{Introduction}
The theory of convex functions has significant applications in optimization and variational inequality problems. Many scientists and researchers have explored it broadly in finite as well as infinite dimensional linear spaces, for the details see[1-7]. Because of its vast applications many generalizations have been proposed. Rapcsak [8] and Udriste [9] presented an innovative generalization, in which the line segment is replaced by the geodesic and Euclidean space is replaced by the Riemannian manifolds. 

\vspace{.2cm}
\noindent
Motivated by the fact that the model of convex function unveils all its applications and outcomes only when it is developed on Riemannian manifolds, we define strong geodesic convex function and  strongly monotone vector field of order $m$ on Riemannian manifolds.

\vspace{.2cm}
Some variational inequality problems on Euclidean spaces can not be solved using classical technique but can be solved on manifolds. Nemeth [3] introduced the concept of variational inequality problem (VIP) on Hadamard manifolds while Li et. al. [10] discussed its existence and uniqueness on Riemannian manifolds. 

\vspace{.2cm}
Motivated by the above mentioned research work, we introduce  variational inequality problem and establish a relation between its solution and the strict minimizers of order $m$ for a multiobjective programming problem.\\

\noindent
\section{Preliminaries}

Let $(M,g)$ be a complete $n$-dimensional Riemannian manifold with Riemannian connection $\nabla$. Let $T_xM$ be the tangent space at $x \in M$ and $\langle.,.\rangle_x$ denotes the scalar product on $T_xM$ with the associated norm denoted by $\| . \|_x$. A vector field $V$ on $M$ is a map from $M$ to $TM$ which associates with each pint $x \in M$ a vector $V(x) \in T_xM$. Let $\gamma_{xy}:[0,1] \rightarrow M$ be a geodesic joining the points $x,y \in M$ such that $\gamma_{xy}(0)=x$ and $\gamma_{xy}(1)=y$. For the basic definitions and concepts of Riemannian geometry one can see ( [11], [12]). 

\vspace{.2cm}
\noindent
Rapcsak [ 8] defined geodesic convexity as follows.

\vspace{.2cm}
\noindent
{\bf Definition 2.1.} [8] A set $A\subseteq M$ is called geodesic convex if a geodesic joining any two points $x,y \in A$ belongs to $A$. 

\vspace{.2cm}
\noindent
{\bf Definition 2.2.} [8] A real valued function $f: A\rightarrow R$  is called geodesic convex if 
\begin{eqnarray*}
f(\gamma_{xy}(t))\leq (1-t)f(x)+tf(y)
\end{eqnarray*}

\noindent
for every $x,y \in A$ and $t \in [0,1]$

\vspace{.2cm}
\noindent
Udriste [9] has defined totally convex set as follows.

\vspace{.2cm}
\noindent
{\bf Definition 2.3.}  [9] A set $A\subseteq M$ is called totally convex if $A$ contains every geodesic $\gamma_{xy}$ of $A$ whose end points $x$ and $y$ are in $A$.

\noindent
\section{Strong geodesic convex functions of order m}

Lin et al. [13 ] extended the concept of convexity to strong convexity of order $m$ on $R^n$ as follows.

\vspace{.2cm}
\noindent
{\bf Definition 3.1.} [13] Let $X$ be a convex seubset of $R^n$. A function $f: X\rightarrow R$ is said to be strongly convex of order $m$ if there exists a constant $c>0$ such that  
\begin{eqnarray*}
f(tx+(1-t)y)\leq tf(x)+(1-t)tf(y)-ct(1-t)\parallel x-y \parallel^m
\end{eqnarray*}
\noindent
for any $x,y \in X$ and $t \in [0,1].$

\vspace{.2cm}
\noindent
Motivated by Lin et al. [13], we introduce the concept of strong geodesic convex function of order $m$.\\

\noindent
{\bf Definition 3.2.}  Suppose $A\subseteq M$ is a geodesic convex set of $M$. A function $f: M\rightarrow R$ is said to be strongly geodesic convex of order $m>0$ on $A$ if there exists a constant $c>0$ such that 
\begin{eqnarray*}
f(\gamma_{xy}(t))\leq (1-t)f(x)+tf(y)-ct(1-t)\parallel \dot {\gamma}_{xy}(t) \parallel^m
\end{eqnarray*}

\noindent
for every $x,y \in A$ and $t \in [0,1]$\\

\noindent
{\bf Remark.} Let $c=0$. Then the above definition becomes the definition of  geodesic convex defined by Rapcsak [8].\\

\noindent
{\bf Theorem 3.1.} Suppose $A\subseteq M$ is a geodesic convex set and $f: M\rightarrow R$ be continuously defferentiable on $A$. Then, $f$ is strongly geodesic convex of order $m$ on $A$ if and only if there exists a constant $c>0$, such that,
$$f(y) \geq f(x) + \dot {\gamma}_{xy}(f)(x) + c \parallel \dot {\gamma}_{xy}(t) \parallel^m,~~~\forall~x,y \in A \eqno(1)$$

\noindent
{\bf Proof.} From the definition of strongly geodesic convex, we have
$$f(\gamma_{xy}(t))\leq (1-t)f(x)+tf(y)-ct(1-t)\parallel \dot {\gamma}_{xy}(t) \parallel^m,\forall~x,y \in A,~~t \in (0,1]$$
\noindent
or
$$f(x)+\frac{f(\gamma_{xy}(t))-f(x)}{t} \leq f(y)-c(1-t)\parallel \dot {\gamma}_{xy}(t) \parallel^m$$
\noindent
Taking limit $t\rightarrow 0$, we get 
$$f(x)+ \dot {\gamma}_{xy}(0)(f) \leq f(y)-c\parallel \dot {\gamma}_{xy}(t) \parallel^m$$
\noindent
or
$$f(y) \geq  f(x)+ \dot{\gamma}_{xy}(f)(x)+ c\parallel \dot {\gamma}_{xy}(t) \parallel^m \eqno(2)$$

\noindent
Conversely, let the given condition holds true for some $c>0$. Changing $y$ with $x$, we get
$$f(x) \geq  f(y)+ \overline {\dot{\gamma}}_{yx}(f)(y)+ c\parallel \overline {\dot {\gamma}}_{yx}(t) \parallel^m \eqno(3)$$

\noindent
where $\overline {\dot {\gamma}}_{yx}(t)=\gamma_{xy}(1-t),~t\in [0,1]$ is a geodesic joining $y$ with $x$. After fixing $t$  we get the point $\gamma_{xy}(t)$. Let $\gamma_{xy}(u)$, $u \in [t,1]$  be the  restriction for the geodesic arc that joins $\gamma_{xy}(t)$ and $y$. 

\noindent
Setting $u=t+s(1-t),~s \in [0,1],$ we obtainn the reparametrization
$$\alpha(s)=\gamma_{xy}(u(s))=\gamma_{xy}(t+s(1-t)),~~s \in [0,1],$$

\noindent
where $\alpha(0)=\gamma_{xy}(t),~~\frac{d\alpha (0)}{ds}=(1-t)\frac{d \gamma_{xy}(t)}{dt}$

\noindent
Similarly, the restriction $\overline {\dot {\gamma}}_{yx}(u)=\gamma_{xy}(1-u),~u\in [1-t,1]$ is a geodesic joining $\gamma_{xy}(t)$ with $x$. Setting $u=(1-t)+st,~s \in [0,1]$, we find the reparametrization
$$\beta(s)=\overline {\gamma}_{yx}(1-t+st)=\gamma_{xy}(t-st),~s \in [0,1],$$

\noindent
where $\beta(0)=\gamma_{xy}(t),~~\frac{d\beta(0)}{ds}=-t \frac{d \gamma_{xy}(t)}{dt}.$ \newline

\vspace{.2mm}
\noindent
On replacing $x$ with $\gamma_{xy}(t)$ in $(2)$ and $\dot {\gamma}_{xy}(0)$ by $\frac{d\alpha (0)}{ds}$, we get
$$f(y) \geq f(\gamma_{xy}(t)) + (1-t) \frac{d \gamma_{xy}(f)}{dt} (\gamma_{xy}(t)) + c \parallel \dot {\gamma}_{xy}(t) \parallel^m \eqno(4)$$
\noindent
Analogously, replacing $y$ with $\gamma_{xy}(t)$ and $\overline {\dot {\gamma}}_{yx}(0)$ by $\frac{d\beta(0)}{ds}$ in $(3)$, we get
$$f(x) \geq f(\gamma_{xy}(t)) -t \frac{d \gamma_{xy}(f)}{dt} (\gamma_{xy}(t)) + c \parallel \overline  {\dot {\gamma}}_{xy}(t) \parallel^m \eqno(5)$$

\noindent
Multiplying $(4)$ by $t$, $(5)$ by $(1-t)$ and then adding, we get
$$f(\gamma_{xy}(t))\leq (1-t) f(x) +tf(y) - ct \parallel \dot {\gamma}_{xy}(t) \parallel^m -c(1-t) \parallel  {\dot {\gamma}}_{xy}(1-t) \parallel^m$$
$$\leq (1-t) f(x) +tf(y) - c't(1-t)\parallel \dot {\gamma}_{xy}(t) \parallel^m~~~~~~~~~~~$$

\noindent
Which shows that $f$ is strongly geodesic convex of order $m$.\\

\noindent
Nemeth [14]  defined monotone vector fields on Riemannian manifolds as follows.

\vspace{.2cm}
\noindent
{\bf Definition 3.3.} Let $M$ be a Riemannian manifold and $V$ be a vecor field on $M$. $V$  is called monotone on $M$ if for every $x,y \in M$ 
$$ \langle V(y), \dot {\gamma}_{xy}(0) \rangle  \leq \langle V(x), \dot {\gamma}_{xy}(1) \rangle$$

\noindent
where $\dot {\gamma}$ denotes the tangent vector of $\gamma$  with respect to the arc length.

\vspace{.2cm}
\noindent
We define strong monotone vector field of order $m$ and establish a relation with strong geodesic convex function of order $m$.\\

\noindent
{\bf Definition 3.4.} Suppose $A \subset M$  is a geodesic convex set. A vector field $V$ on $A$ is called strongly monotone of order $m$ if there exists a constant $\beta > 0$ such that
$$ \langle V(y), \dot {\gamma}_{xy}(1) \rangle - \langle V(x), \dot {\gamma}_{xy}(0) \rangle \geq \beta \parallel \dot {\gamma}_{xy}(t) \parallel^m,~~~~~~~\forall~x,y \in A.\eqno(6)$$

\noindent
{\bf Theorem 3.2.} Suppose $A \subset M$  is a geodesic convex set and $f: M\rightarrow R$ be continuously defferentiable on $A$. Then, $f$ is strongly geodesic convex of order $m$ on $A$ iff $\dot {\gamma}_{xy}(f)$ is strongly monotone of order $m$ on $A$.\\

\noindent
{\bf Proof.} Let $f$ be strongly geodesic convex of order $m$ on $A$. By Theorem 3.1, there exists a constant $c>0$ such that $(1)$ holds. Then, for any $x,y \in A$, we have
$$f(y)-f(x)\geq  \dot {\gamma}_{xy}(f)(x) + c \parallel \dot {\gamma}_{xy}(t) \parallel^m,$$
$$f(x)-f(y)\geq  \dot {\gamma}_{yx}(f)(y) + c \parallel \dot {\gamma}_{xy}(t) \parallel^m,$$
\noindent
On adding, we get
$$\dot {\gamma}_{xy}(f)(x)+\dot {\gamma}_{yx}(f)(y)+2c\parallel \dot {\gamma}_{xy}(t) \parallel^m \leq 0$$
\noindent
or
$$\dot {\gamma}_{yx}(f)(x)-\dot {\gamma}_{yx}(f)(y) \geq 2c\parallel \dot {\gamma}_{xy}(t) \parallel^m $$

\noindent
Which shows that if part is true.

\vspace{.2cm}
\noindent
Conversely, let $(6)$ holds true and $V=\nabla f$. Set $t_i=\frac{i}{m+1},~i=0,1,2,....m+1.$

\vspace{.2mm}
\noindent
By the mean value theorem, $\exists~\xi \in (t_i, t_{i+1}),~0 \leq i \leq m$

$$f({\gamma}_{xy}(0)+t_{i+1}({\gamma}_{xy}(1)-{\gamma}_{xy}(0))-f({\gamma}_{xy}(0)+t_i({\gamma}_{xy}(1)-{\gamma}_{xy}(0))$$ $$= (t_{i+1}-t_i)({\gamma}_{xy}(1)-{\gamma}_{xy}(0) )^T \nabla f(x+\xi_i ({\gamma}_{xy}(1)-{\gamma}_{xy}(0)))$$
\noindent
It follows from $(6)$,

$$\begin{array}{lll}
f(y)-f(x) &=&\Sigma_{i=0}^m [f({\gamma}_{xy}(0)+t_{i+1}({\gamma}_{xy}(1)-{\gamma}_{xy}(0))-f({\gamma}_{xy}(0)+t_i({\gamma}_{xy}(1)-{\gamma}_{xy}(0)))] \\
\\
&=& \Sigma_{i=0}^m (t_{i+1}-t_i) ({\gamma}_{xy}(1)-{\gamma}_{xy}(0) )^T [\nabla f(x+\xi_i ({\gamma}_{xy}(1)-{\gamma}_{xy}(0))-\nabla f(x))]
\\
&+&({\gamma}_{xy}(1)-{\gamma}_{xy}(0) )^T \nabla f(x))\\
\\
&\geq& \beta \parallel \dot {\gamma}_{xy}(t) \parallel^m \Sigma_{i=0}^m \xi_i^{m-1} (t_{i+1}-t_i) + ({\gamma}_{xy}(1)-{\gamma}_{xy}(0) )^T \nabla f(x)
\end {array}$$

\noindent
Taking limit $m\rightarrow\infty$, we get
$$f(y)-f(x) \geq \frac{\beta}{m} \parallel \dot {\gamma}_{xy}(t) \parallel^m + ({\gamma}_{xy}(1)-{\gamma}_{xy}(0) )^T \nabla f(x).$$

Using theorem 3.1, the result follows.\\

\noindent
{\bf Definition 3.5.} Suppose $A \subset M$ is a geodesic convex set. A vector field $V$ on $A$ is called strongly pseudomonotone of order $m$ if
$$ \langle V(x), \dot {\gamma}_{xy}(0) \rangle+ \beta \parallel \dot {\gamma}_{xy}(t)\rangle \parallel^m \geq 0 \Rightarrow  \langle V(y), \dot {\gamma}_{xy}(1) \rangle \geq 0 ~~~~~~~\forall~x,y \in A. \eqno(7)$$

\noindent
{\bf Proposition 3.1.} Every strongly monotone vector field of order $m$ is strongly pseudomonotone of order $m$.\\

\noindent
{\bf Proof.} Let $V$ on $A$ be strongly monotone of order $m$, then 
$$ \langle V(y), \dot {\gamma}_{xy}(1) \rangle - \langle V(x), \dot {\gamma}_{xy}(0) \rangle \geq \beta \parallel \dot {\gamma}_{xy}(t) \parallel^m$$

\noindent
or
$$ \langle V(y), \dot {\gamma}_{xy}(0)\rangle +\beta \parallel \dot {\gamma}_{xy}(t) \parallel^m \leq \langle V(x), \dot {\gamma}_{xy}(1) \rangle $$

\noindent
Let $\langle V(y), \dot {\gamma}_{xy}(0)\rangle+\beta \parallel \dot {\gamma}_{xy}(t) \parallel^m  \geq 0$, then
$$\langle V(x), \dot {\gamma}_{xy}(1) \rangle \geq 0.$$ 
\noindent
Hence, $V$ is strongly pseudomonotone of order $m$.


\noindent
\section{Variational Inequality Problem}

\noindent
Let $M$ be a complete Riemannian manifold and $A \subseteq M$ be a non empty set of $M$. Let $\Gamma_{x,y}^A$ denotes the collection of all geodesics from $x$ to $y$  such that $\gamma_{xy}\in A$. Suppose that $T=(T_1,T_2,....T_k)$, where $T_i: A\rightarrow 2^{TM}$ be a set-valued vector field on $A$. The variational inequality problem is to find $\overline x \in A$ and $v \in T_{\overline x} M$ such that 
$$\langle v, \dot {\gamma}_{\overline x x}(0) \rangle \not< 0~~~~{\mbox for~all} x \in A,$$ 

\noindent
where $\langle v, \dot {\gamma}_{\overline x x}(0) \rangle = (\langle v_1, \dot {\gamma}_{\overline x x}(0) \rangle, \langle v_2, \dot {\gamma}_{\overline x x}(0) \rangle....,\langle v_k, \dot {\gamma}_{\overline x x}(0) \rangle)$, $v_i \in T_{i_{\overline x}}M$, $i=1,2,....k$.

\noindent
The multiobjective optimization problem (MOP) is to find a strict minimizer of order $m$ for 
$$ {\mbox minimize}~f(x)=(f_1(x), f_2(x),....f_k(x)),~~~~x \in A$$

\noindent
{\bf Theorem 4.1.} Let $f_i,~i=1,2....,k$ be strongly convex of order m on $A$. Then $\overline x \in A$ is the solution of VIP with $T_{i_{\overline x}}M=\nabla f_i(\overline x),~i=1,2,....k$, iff $\overline x$ is a strict minimizer of order $m$ for the MOP.\\

\noindent
{\bf Proof.} Suppose $\overline x$ is the solution of VIP but is not a strict minimizer of order $m$ for MOP. Then for $c>0$, there exists some $x^* \in A$, such that
$$f(x^*) < f(\overline x)+c \parallel \dot {\gamma}_{\overline x x^*}(0) \parallel^m$$
or
\noindent
$$  f_i(x^*) < f_i(\overline x)+c_i \parallel \dot {\gamma}_{\overline x x^*}(0) \parallel^m,~i=1,2,....k.$$

\noindent
Since $f_i,~i=1,2,....k$, are strongly geodesic convex of order $m$ on $A$, the above inequality implies
$$\langle v_i, \dot {\gamma}_{\overline x x^*}(0) \rangle <0,~~\forall~v_i \in T_{i_{\overline x}}M=\nabla f_i(\overline x),~i=1,2,....k,$$
\noindent
that is
$$\langle v, \dot {\gamma}_{\overline x x^*}(0) \rangle < 0 ~~\forall~v \in T_{\overline x}M,~x^* \in A,$$
\noindent
which contradicts the assumption that $\overline x$ is the solution of the VIP.

\parindent=8mm
Conversely, let $\overline x$ be a strict minimizer of order $m$ for (MOP) but is not a solution of (VIP). Therefore, there exists an $x^* \in A$ such that 
$$\langle v_i, \dot {\gamma}_{\overline x x^*}(0) \rangle <0,~~\forall~v_i \in T_{i_{\overline x}}M=\nabla f_i(\overline x),~i=1,2,....k.$$
\noindent
Using the definition of strongly geodesic convexity of order $m$ for $f_i, i=1,2,....k,$ we get 
$$ f_i(x^*)-f_i(\overline x) < c \parallel \dot {\gamma}_{\overline x x^*}(t) \parallel^m,$$
\noindent
which contradicts the strict minimizer condition. Consequently, $\overline x$ is not a strict minimizer of order $m$ for (MOP) and hence $\overline x$ is a solution of (VIP).\\

\vspace{.4cm}
\parindent=8mm
{{ \bf REFERENCES}
\begin{enumerate}

\item Pini, R.: Convexity along curves and invexity. Optimization 29, 301-309 (1994).
\item Mititelu, S.: Generalized invexity and vector optimization on differentiable manifolds. Diff. Geom. Dyn. Syst. 3, 21-31 (2001).
\item Nemeth, S.Z.: Variational Inequalities on Hadamard Manifolds. Nonlinear Anal. 52: 1491-1498 (2003).
\item  Iqbal, A., Ali, S., Ahmad,  I.,  On geodesic E-convex sets, geodesic E-convex functions and E-epigraphs. J Optim Theory Appl 155(1): 239-251 (2012).
\item  Iqbal, A., Ahmad,  I., Ali, S., Some properties of geodesic semi-E-convex functions. Nonlinear Anal. 74, 6805-6813 (2011).
\item Iqbal, A.: On  $\varphi_h$-preinvex functions. Comm. Appl. Anal. 20: 175-185 (2016).
\item  Pouryayevali,  M.R., Barani, A.,  Invariant monotone vector fields on Riemannian manifolds, Nonlinear Anal. 70 (5), 1850-1861 (2009).

\item Rapcsak, T.: Smooth Nonlinear Optimization in $\mathbb{R}^n$, Kluwer Academic Publishers (1997).
\item  Udriste, C.: Convex Functions and Optimization Methods on Riemannian Manifolds, Kluwer Academic Publishers (1994).
\item Li, S.L., Li, C., Liou, Y.C., Yao, J.C.: Existence of solutions for variational inequalities on Riemannian manifolds. Nonlinear Anal. 71, 5695-5706 (2009).
\item do Carmo, M.P.: Riemannian Geometry, Bikhauser, Boston, 1992.
\item Spivak, M.: Calculus on Manifolds: A Modern Aproach to Classical Theorems of Advanced Calculus, W.A. Benjamin, Inc, New York, Amsterdam, 1965.
\item Lin, G.H., Fukushima, M.: Some exact penalty results for nonlinear programs and mathematical programs with equilibrium constraints. J Optim Theory Appl. 118 (1), 67-80 (2003).
\item Nemeth, S.Z.: Monotone vector fields, Publ. Mat. 54, 437-449 (1999).
\end{enumerate}

\end{document}